\documentclass[a4paper,10pt]{article}
\usepackage[top=2.5cm, left=2cm, right=2cm, bottom=2.5cm]{geometry}
\usepackage{float}
\floatstyle{ruled}
\usepackage{listings}
\usepackage{amsfonts, booktabs}
\usepackage{verbatim,amsmath,amssymb,graphicx}
\usepackage{url}

\newfloat{example}{!ht}{loa}[section]
\floatname{example}{}

\newcommand\fining{{\sf FinInG}}
\newcommand\calS{\mathcal{S}}
\newcommand\cE{\mathcal{E}}
\newcommand\cC{\mathcal{C}}
\newcommand\pg{\mathrm{PG}}
\newcommand\q{\mathrm{Q}}
\newcommand\w{\mathrm{W}}
\newcommand\h{\mathrm{H}}
\newcommand\gf{\mathrm{GF}}
\newcommand\scalars{\mathrm{Sc}}
\newcommand\gammal{\Gamma\mathrm{L}}
\newcommand\pgammal{\mathrm{P}\Gamma\mathrm{L}}
\newcommand\agammal{\mathrm{A}\Gamma\mathrm{L}}
\newcommand\pgl{\mathrm{PGL}}
\newcommand\SL{\mathrm{SL}}
\DeclareMathOperator\I{\mathrm{I}}

\usepackage{xcolor}
\title{\fining{}: a package for Finite Incidence Geometry}

\author{John Bamberg, Anton Betten, Philippe Cara, Jan De Beule, \\Max Neunh\"offer, and Michel Lavrauw}

\begin{document}

\maketitle

\begin{abstract}
\fining{} is a package for computation in Finite Incidence Geometry. It provides users with the basic tools to work 
in various areas of finite geometry from the realms of projective spaces to the flat lands of generalised polygons. 
The algebraic power of GAP is exploited, particularly in its facility with matrix and permutation groups.
\end{abstract}



\tableofcontents

\section{Finite incidence geometry on the computer}	

\fining{} \cite{fining} is a new share package for GAP \cite{GAP4}. It is designed to provide an environment to explore various finite geometrical structures,
such as projective spaces, affine spaces, classical polar spaces, coset geometries and their diagrams and generalised polygons; 
allowing an integrated use of the automorphism groups and exploring geometrical relations between different geometries via incidence 
preserving morphisms. The use of automorphism groups is connected with the existing primary functionality of GAP as a computer algebra package chiefly concerned
with computing with groups. The algebraic power of GAP is especially exploited when finite fields are involved, and when considering affine and projective algebraic
varieties. Our package interfaces with the packages {\sf GRAPE} \cite{grape} and {\sf DESIGN} \cite{design} which is especially useful when working with
generalised polygons and purely combinatorial incidence structures. The representation of an automorphism group of a projective space
relies on an implementation of collineations of a projective space as semilinear maps of the underlying vector space, and this is done in conjunction with the 
packages {\sf GenSS} and {\sf Orb} \cite{genss,orb}  that provide sophisticated algorithms for enumeration of orbits, stabilisers, and stabiliser chains. 
Finally the package {\sf Forms} \cite{forms}
enables us to provide the user with the possibility to construct classical polar spaces from a user-defined sesquilinear or quadratic form,
computing the necessary coordinate transformations on the fly. This approach results in a user-friendly and intuitive package
providing generic functionality for exploring the most commonly studied finite incidence structures and their automorphism groups.

Other good systems are available for incidence geometry but often focus on the efficient handling of particular cases, 
rather than on the generic approach which is central in \fining.
The system {\sf MAGMA} \cite{magma}, for example, is renown as a system for computational group theory. It contains
an implementation of projective planes, some substructures, coset geometries and related functionality. 
The computation of the automorphism group of a projective plane is done very efficiently, however, 
it is always returned as a permutation group on the points of the plane. The creation of a coset geometry is based on a given permutation group, 
which is very efficient, but which makes it more difficult for the user to use, for example, matrix groups occurring from (say) the 
automorphism group of a classical polar space. For the student 
starting with finite incidence geometry, or the (experienced) researcher interested
in exploring a particular incidence structure using a computer, this approach might be a hurdle.

This explains the basic philosophy of \fining{}: to construct an incidence structure (together with an automorphism group) 
in its natural representation, where the user decides what \emph{natural} actually means! For example, a conic
with a particular equation can easily be 
constructed, or a coset geometry arising from a matrix group or finitely presented group instead of a permutation group.
Once objects in a certain category have been constructed, their behaviour is independent 
of their representation. Furthermore, many well studied geometries are considered in their standard representation, 
together with straight-forward methods for switching from one representation to another.
\begin{example}
{\footnotesize
\begin{lstlisting}[language=GAP]
gap> map := KleinCorrespondence(7);
<geometry morphism from <lines of ProjectiveSpace(3, 7)> to 
<points of Q+(5, 7): x_1*x_6+x_2*x_5+x_3*x_4=0>>
gap> r := Range(map);
<points of Q+(5, 7): x_1*x_6+x_2*x_5+x_3*x_4=0>
gap> hyp := AmbientGeometry(r);
Q+(5, 7): x_1*x_6+x_2*x_5+x_3*x_4=0
gap> klein := HyperbolicQuadric(5,7);
Q+(5, 7)
gap> EquationForPolarSpace(klein);
x_1*x_2+x_3*x_4+x_5*x_6
gap> map := KleinCorrespondence(klein);
<geometry morphism from <lines of ProjectiveSpace(3, 7)> to <points of Q+(5, 7)>>
gap> morphism := IsomorphismPolarSpaces(hyp,klein);
<geometry morphism from <Elements of Q+(5, 7): x_1*x_6+x_2*x_5+x_3*x_4=0> to 
<Elements of Q+(5, 7)>>
\end{lstlisting}
\caption{Klein Correspondence}\label{ex:klein}%
}
\end{example}
\begin{example}
{\footnotesize
\begin{lstlisting}[language=GAP]
gap> L := SimpleLieAlgebra("E",6,Rationals);
<Lie algebra of dimension 78 over Rationals>
gap> rs := RootSystem(L);
<root system of rank 6>
gap> w := WeylGroup(rs);
<matrix group with 6 generators>
gap> gens := GeneratorsOfGroup(w);;
gap> parabolics := List(gens, x -> Group(Difference(gens, [x])));
[ <matrix group with 5 generators>, <matrix group with 5 generators>, 
  <matrix group with 5 generators>, <matrix group with 5 generators>, 
  <matrix group with 5 generators>, <matrix group with 5 generators> ]
gap> g := Group(gens);
<matrix group with 6 generators>
gap> cg := CosetGeometry(g,parabolics);
CosetGeometry( Group(
[ [ [ -1, 0, 1, 0, 0, 0 ], [ 0, 1, 0, 0, 0, 0 ], [ 0, 0, 1, 0, 0, 0 ], 
      [ 0, 0, 0, 1, 0, 0 ], [ 0, 0, 0, 0, 1, 0 ], [ 0, 0, 0, 0, 0, 1 ] ], 
  [ [ 1, 0, 0, 0, 0, 0 ], [ 0, -1, 0, 1, 0, 0 ], [ 0, 0, 1, 0, 0, 0 ], 
      [ 0, 0, 0, 1, 0, 0 ], [ 0, 0, 0, 0, 1, 0 ], [ 0, 0, 0, 0, 0, 1 ] ], 
  [ [ 1, 0, 0, 0, 0, 0 ], [ 0, 1, 0, 0, 0, 0 ], [ 1, 0, -1, 1, 0, 0 ], 
      [ 0, 0, 0, 1, 0, 0 ], [ 0, 0, 0, 0, 1, 0 ], [ 0, 0, 0, 0, 0, 1 ] ], 
  [ [ 1, 0, 0, 0, 0, 0 ], [ 0, 1, 0, 0, 0, 0 ], [ 0, 0, 1, 0, 0, 0 ], 
      [ 0, 1, 1, -1, 1, 0 ], [ 0, 0, 0, 0, 1, 0 ], [ 0, 0, 0, 0, 0, 1 ] ], 
  [ [ 1, 0, 0, 0, 0, 0 ], [ 0, 1, 0, 0, 0, 0 ], [ 0, 0, 1, 0, 0, 0 ], 
      [ 0, 0, 0, 1, 0, 0 ], [ 0, 0, 0, 1, -1, 1 ], [ 0, 0, 0, 0, 0, 1 ] ], 
  [ [ 1, 0, 0, 0, 0, 0 ], [ 0, 1, 0, 0, 0, 0 ], [ 0, 0, 1, 0, 0, 0 ], 
      [ 0, 0, 0, 1, 0, 0 ], [ 0, 0, 0, 0, 1, 0 ], [ 0, 0, 0, 0, 1, -1 ] ] ]) )
\end{lstlisting}
\caption{An apartment in a building of type $E_6$}\label{ex:e6}%
}
\end{example}
Example~\ref{ex:klein} illustrates the Klein correspondence (between the lines of a $3$-dimensional projective 
space and the points of a hyperbolic quadric in $5$ dimensions). The example shows the default 
behaviour and the possibility to employ a user-defined form for the Klein quadric.

Example~\ref{ex:e6}, illustrates the construction of an apartment in a building of type $\mathrm{E}_6$, based on existing 
GAP functionality of root systems and Lie algebras.

\section{Fundamentals of \fining{}}

First of all, a set of new data structures, which are type in GAP by the so-called {\em categories}, is provided in \fining{},
summarized as follows.
\begin{itemize}
\item \fining{} Incidence geometries have their category, together with elementary operations.
\item Finite projective spaces, classical finite polar spaces, and some classical generalized polygons are incidence geometries
whose elements are represented by sub spaces of a vector space over a finite field. They are called Lie geometries in \fining{}.
More particular operations are applicable on such geometries. As explained in the introduction, the user is free to choose a 
representation for a Lie geometry.
\item Generalised polygons are particular examples of point-line incidence geometries. Many examples in many different representations
exist. \fining{} provides a generic category, together with some particular operations, allowing the user to construct in a very natural way
their desired geometry.
\item Coset geometries can be constructed from any group in their own category. 
\item Finally affine spaces can also be constructed as objects in their dedicated category. 
\end{itemize}

The automorphism groups of Lie geometries are naturally represented by semi-linear maps on the underlying vector space. \fining{}
provides category for such maps acting on the projective space. Groups generated by such maps can be constructed, functions representing 
the action on elements of a Lie geometry are provided, and actions can be computed, including nice monomorphisms. This part of 
\fining{} is the direct link between the automorphism groups of Lie geometries and the available GAP operations for permutation groups.

Secondly fining provides efficient operations to compute orbits and (setwise) stabilisers of projective semilinear groups. These
operations are based on the use of the packages orb, genss and cvec. Also from the algorithmic point of view, \fining{} provides
efficient methods to enumerate the the elements of finite classical polar spaces. Geometry morphisms can be constructed. This 
is especially powerfull combined with field reduction.

\section{Incidence geometries}

We follow \cite{BC13} for the definitions of incidence structure and incidence geometry. An {\em incidence structure}  $\calS$
consists of a set of elements, a reflexive and symmetric relation $\I$ on the elements and a type function from the set of 
elements to an index set (i.e., every element has a ``type''). It satisfies the following axiom: 
\begin{quote}
(i) no two different elements of the same type are incident. 
\end{quote}
If we forget the type function, an incidence structure is just a multipartite graph where adjacency is incidence.
Note that every element is incident with itself, so there is a loop on each vertex.
The type function turns the graph into a coloured graph. 
Assume that $\calS$ is an incidence structure with type function $t$ with range $T = \{1,\ldots,n\}$. Then we say that $\calS$ has rank $n$.
A {\em flag} is a set $F$ of elements of $\calS$ such that every two elements contained in $F$ are incident. Necessarily, 
$F$ cannot contain two elements of the same type, so $F$ contains at most $n$ elements. We call a flag {\em maximal} when
no element of $\calS$ is incident with all elements of $F$, unless it is
already contained in $F$.

The term geometry, or incidence geometry, is interpreted broadly in this package. Particularly, an {\em incidence geometry} is an 
incidence structure satisfying the following axiom: 
\begin{quote}
(ii) every maximal flag contains an element of each type. 
\end{quote}
In graph terminology, this means that every maximal clique contains an element of each colour.

Let $\calS$ be an incidence structure of rank $n$, and let $e$ be an element of type $i \in T$. Let $j \in T \setminus \{i\}$. 
Then the {\em $j$-shadow of $e$} is the set of elements of type $j$ incident with $e$. Let $F$ be a flag of $\calS$ then
the {\em $j$-shadow of $F$} is the set of elements of type $j$ incident with all elements in $F$. The {\em residue of $F$}
is the set of all elements of $\calS$ not in $F$ but incident with all elements of $F$, together with the induced geometrical structure. 
Hence, the residue of a maximal flag of an incidence geometry 
is empty. All geometries that can be constructed in \fining{} are incidence structures. This terminology is consistently used, as illustrated
in Example~\ref{ex:cosetelements} for an ``abstract'' coset geometry and Example~\ref{ex:elements_proj} for a projective
space.
\begin{example}
{\footnotesize
\begin{lstlisting}[language=GAP]
gap> involset:=[(1,2), (1,3), (3,4), (4,5), (4,6), (6,9), (6,10),
(7,10), (8,10)];;
gap> g:=Group(involset);;
gap> parabolics := List(involset, i -> Subgroup(g, Difference(involset, [i]) ) );;
gap> img:=CosetGeometry(g, parabolics);
CosetGeometry( Group( [ (1,2), (1,3), (3,4), (4,5), (4,6), (6,9),
( 6,10), ( 7,10), ( 8,10) ] ) )
gap> elements_type_2 := ElementsOfIncidenceStructure(img,2); 
<elements of type 2 of CosetGeometry( Group( 
[ (1,2), (1,3), (3,4), (4,5), (4,6), (6,9), ( 6,10), ( 7,10), ( 8,10) ]
) )>
gap> el2 := Random(elements_type_2);
<element of type 2 of CosetGeometry( Group( 
[ (1,2), (1,3), (3,4), (4,5), (4,6), (6,9), ( 6,10), ( 7,10), ( 8,10) ]
) )>
gap> shad := ShadowOfElement(img,el2,3);
<shadow elements of type 3 in CosetGeometry( Group( 
[ (1,2), (1,3), (3,4), (4,5), (4,6), (6,9), ( 6,10), ( 7,10), ( 8,10) ]
) )>
gap> el3 := Random(shad);
<element of type 3 of CosetGeometry( Group( 
[ (1,2), (1,3), (3,4), (4,5), (4,6), (6,9), ( 6,10), ( 7,10), ( 8,10) ]
) )>
gap> IsIncident(el2,el3);
true
gap> flag := FlagOfIncidenceStructure(img,[el2,el3]);;
gap> shad := ShadowOfFlag(img,flag,4);
<shadow elements of type 4 in CosetGeometry( Group( 
[ (1,2), (1,3), (3,4), (4,5), (4,6), (6,9), ( 6,10), ( 7,10), ( 8,10) ]
) )>
\end{lstlisting}
\caption{Elements and shadows in a coset geometry}\label{ex:cosetelements}%
}
\end{example}

\begin{example}
\begin{lstlisting}[language=GAP]
gap> pg := PG(5,7);
ProjectiveSpace(5, 7)
gap> Rank(pg);
5
gap> line := Random(Lines(pg));
<a line in ProjectiveSpace(5, 7)>
gap> shadowplanes := ShadowOfElement(pg,line,3); #planes are type 3
<shadow planes in ProjectiveSpace(5, 7)>
gap> plane := Random(shadowplanes);
<a plane in ProjectiveSpace(5, 7)>
gap> planepoints := ShadowOfElement(pg,plane,1); #planes are type 1
<shadow points in ProjectiveSpace(5, 7)>
gap> p := Random(planepoints);
<a point in ProjectiveSpace(5, 7)>
gap> flag := FlagOfIncidenceStructure(pg,[p,plane]);
<a flag of ProjectiveSpace(5, 7)>
gap> shadowlines := ShadowOfFlag(pg,flag,2);
<shadow lines in ProjectiveSpace(5, 7)>
\end{lstlisting}
\caption{Elements and shadows in a projective space}\label{ex:elements_proj}
\end{example}

\section{Lie geometries}

Example~\ref{ex:elements_proj} brings us to the Lie geometries. A {\em Lie geometry}, is an incidence geometry 
whose automorphism group is a (subgroup of) a group of Lie type in its natural representation. One common fact about Lie geometries is 
that their elements are represented by subspaces of a vector space. All the collineations are induced by semilinear 
maps of the underlying vector space. So a finite Lie geometry is a geometry for which the elements are subspaces  
of some projective space over a finite field, and with symmetrised
set-theoretic containment as incidence. Hence, a Lie geometry is 
naturally embedded in a projective space. These principles are used in \fining{} to implement Lie geometries. 
Currently the following Lie geometries are readily implemented in \fining{}:
\begin{itemize}
\item projective spaces,
\item classical polar spaces, and,
\item the so-called classical generalised hexagons and octagons.
\end{itemize}

Elements of Lie geometries are also called {\em subspaces}. Subspaces have a (projective) dimension, and subspaces of projective
dimension $0,1,2,3$ are called points, lines, planes, solids, respectively. In \fining{} the default type set is $\{1,\ldots,n\}$ 
and the type of an element is its projective dimension plus one. Some handy shortcuts for the operation 
\verb|ShadowOfElements| are implemented, so that given a point $p$, \verb|Lines(p)| is equivalent with \verb|ShadowOfElement(pg,p,2)|. 
Basic operations such as \verb|Meet| and \verb|Span|, illustrated in Example~\ref{ex:proj_els}, rely on operations to deal with subspaces 
of the underlying vector space of the Lie geometry.

\begin{example}
{\footnotesize
\begin{lstlisting}[language=GAP]
gap> pg := PG(8,9);
ProjectiveSpace(8, 9)
gap> ElementsOfIncidenceStructure(pg,1);
<points of ProjectiveSpace(8, 9)>
gap> Points(pg);
<points of ProjectiveSpace(8, 9)>
gap> Lines(pg);
<lines of ProjectiveSpace(8, 9)>
gap> Planes(pg);
<planes of ProjectiveSpace(8, 9)>
gap> Solids(pg);
<solids of ProjectiveSpace(8, 9)>
gap> ElementsOfIncidenceStructure(pg,5);
<proj. 4-subspaces of ProjectiveSpace(8, 9)>
gap> plane := Random(Planes(pg));
<a plane in ProjectiveSpace(8, 9)>
gap> Solids(plane);
<shadow solids in ProjectiveSpace(8, 9)>
gap> line := Random(Lines(pg));
<a line in ProjectiveSpace(8, 9)>
gap> Span(line,plane);
<a proj. 4-space in ProjectiveSpace(8, 9)>
gap> Meet(line,plane);
< empty subspace >
gap> hyperplane := Random(Hyperplanes(pg));     
<a proj. 7-space in ProjectiveSpace(8, 9)>
gap> Meet(line,hyperplane);
<a point in ProjectiveSpace(8, 9)>
gap> Span(line,hyperplane);
ProjectiveSpace(8, 9)
\end{lstlisting}
\caption{Lie geometries: elements of projective spaces}\label{ex:proj_els}%
}
\end{example}

A {\em polar space} is a partial linear space satisfying the so-called one-or-all axiom: given a point $P$ and a line $l$,
$P$ not incident with $l$, there is either exactly one line incident with $P$ meeting $l$ in a point, or points of $l$
span with $P$ a line of the polar space. The {\em classical polar spaces} are constructed from vector spaces
equipped with a form as follows. Consider a vector space $V(d+1,q)$ equipped with a non-degenerate sesquilinear form $f$, 
respectively non-singular quadratic form $f$. The polar space $\calS$ associated to the form $f$ is the geometry 
of which the elements are represented by the totally isotropic, respectively totally singular, 
subspaces with respect to the form $f$. The rank of $\calS$ equals the Witt index of the form $f$. So we see that $\calS$ is embedded
in the projective space $\pg(d,q)$. Sesquilinear and quadratic forms
on a finite vector space have been classified. Up to coordinate transformation,
they fall into three main classes.
\begin{enumerate}
\item orthogonal polar spaces, associated with a quadratic form (which in odd characteristic is equivalent with a symmetric bilinear form);
\item symplectic polar spaces, associated with an alternating bilinear form;
\item hermitian polar spaces, associated with a unitary form.
\end{enumerate} 
The orthogonal polar spaces separate into three classes: parabolic quadrics (in even projective dimension), 
and elliptic and hyperbolic quadrics (in odd projective dimension). 

Classical polar spaces can be constructed in several ways in \fining{}. A user-defined sesquilinear or quadratic form can be used, but
\fining{} also provides shortcuts to classical polar spaces defined by a ``standard form''. The following example illustrates this possibility.
When appropriate, the connection with polarities of projective spaces can be explored. Example~\ref{ex:polar} illustrates some basic functionality
for polar spaces as well.

\begin{example}
{\footnotesize
\begin{lstlisting}[language=GAP]
gap> mat := [[0,0,0,1],[0,0,-2,0],[0,2,0,0],[-1,0,0,0]]*Z(5)^0;
[ [ 0*Z(5), 0*Z(5), 0*Z(5), Z(5)^0 ], [ 0*Z(5), 0*Z(5), Z(5)^3, 0*Z(5) ], 
  [ 0*Z(5), Z(5), 0*Z(5), 0*Z(5) ], [ Z(5)^2, 0*Z(5), 0*Z(5), 0*Z(5) ] ]
gap> form := BilinearFormByMatrix(mat,GF(25));
< bilinear form >
gap> ps := PolarSpace(form);
<polar space in ProjectiveSpace(
3,GF(5^2)): x1*y4+Z(5)^3*x2*y3+Z(5)*x3*y2-x4*y1=0 >
gap> Rank(ps);
2
gap> ps;
W(3, 25): x1*y4+Z(5)^3*x2*y3+Z(5)*x3*y2-x4*y1=0
gap> ps := EllipticQuadric(5,7);
Q-(5, 7)
gap> Rank(ps);
2
gap> phi := PolarityOfProjectiveSpace(ps);
<polarity of PG(5, GF(7)) >
gap> ps := ParabolicQuadric(4,4);
Q(4, 4)
gap> phi := PolarityOfProjectiveSpace(ps);
Error, no polarity of the ambient projective space can be associated to <ps>
\end{lstlisting}
\caption{Finite classical polar spaces: elementary operations}\label{ex:polar}%
}
\end{example}
Lie geometries are {\em naturally} embedded in projective spaces. In \fining{} this means that an element of a polar space, for example, is
automatically an element of the ambient projective space. Example~\ref{ex:lie_embedded} demonstrates the philosophy of 
Lie geometries in \fining{}.

\begin{example}
\begin{lstlisting}[language=GAP]
gap> hexagon := SplitCayleyHexagon(3);
H(3)
gap> p := Random(Points(hexagon));
#I  for Split Cayley Hexagon
#I  Computing nice monomorphism...
#I  Found permutation domain...
<a point in H(3)>
gap> lines := Lines(p);
<shadow lines in H(3)>
gap> Size(lines);   
4
gap> polarspace := AmbientPolarSpace(hexagon);
Q(6, 3): -x_1*x_5-x_2*x_6-x_3*x_7+x_4^2=0
gap> p2 := ElementToElement(polarspace,p);
<a point in Q(6, 3): -x_1*x_5-x_2*x_6-x_3*x_7+x_4^2=0>
gap> p2 in hexagon;
true
gap> lines2 := Lines(p2);
<shadow lines in Q(6, 3): -x_1*x_5-x_2*x_6-x_3*x_7+x_4^2=0>
gap> Size(lines2);
40
gap> Collected(List(lines2,x->x in hexagon));
[ [ true, 4 ], [ false, 36 ] ]
gap> pg := AmbientSpace(hexagon);
ProjectiveSpace(6, 3)
gap> p3 := ElementToElement(pg,p);        
<a point in ProjectiveSpace(6, 3)>
gap> p3 in polarspace;
true
gap> p3 in hexagon;
true
gap> lines3 := Lines(p3);
<shadow lines in ProjectiveSpace(6, 3)>
gap> Collected(List(lines3,x->x in hexagon));
[ [ true, 4 ], [ false, 360 ] ]
gap> Collected(List(lines3,x->x in polarspace));
[ [ true, 40 ], [ false, 324 ] ]
\end{lstlisting}
\caption{Lie geometries embedded in projective spaces}\label{ex:lie_embedded}
\end{example}

\section{Collineation groups of Lie geometries}

A {\em collineation} of an incidence geometry $\calS$ is a type preserving map from $\calS$ to itself that preserves the incidence. 
Let $\pg(d,q)$ denote a $d$-dimensional projective space over the finite field $\gf(q)$. As we have seen in the previous section, the elements
of $\pg(d,q)$ are represented by the subspaces of the $d+1$-dimensional vector space $V(d+1,q)$. By the Fundamental Theorem of 
Projective Geometry, each collineation of a projective space $\pg(d,q)$ is induced by a semilinear automorphism of $V(d+1,q)$. 
The group of semilinear automorphisms of $V(d+1,q)$ is denoted by $\gammal(d+1,q)$, and an element of this group is completely
determined by a non-singular matrix and an automorphism of the underlying field. The scalar matrices form a subgroup $\scalars(d+1,q)$
of the group $\gammal(d+1,q)$ that acts trivially on the subspaces of $V(d+1,q)$. The permutation group
induced on the subspaces is the \emph{projective general semilinear group}:
\[ \pgammal(d+1,q) \cong \gammal(d+1,q) / \scalars(d+1,q)\,. \]
The group $\pgl(d+1,q)$ is the subgroup of $\gammal(d+1,q)$ induced by the linear automorphisms of $V(d+1q)$. This group is called the
{\em projectivity group} or {\em homography group} of $\pg(d,q)$. 

The projectivity groups are already implemented in GAP, but only as permutation groups.
without the underlying geometry. 
W provide the implementation of collineation groups of projective spaces as groups 
induced by semilinear transformations of vector spaces. Note that an element of a 
Lie geometry is represented by a subspace of a vector space, which in \fining{} is inturn represented as the row-space of a matrix. 
Let $M$ be an $(m + 1) \times (d+1)$-matrix over $\gf(q)$ of rank $m+1$. The row-space of $M$ represents 
a subspace $\pi$ of $\pg(d,q)$ of projective dimension $m$. Let $\phi$ be a semilinear automorphism of $V(d+1,q)$ represented by the matrix 
$A$ and the field automorphism $\theta$. In \fining{}, the image of $\pi$ under $\phi$ is the projective subspace represented by the matrix $M'$, with
\[ M' := (M A)^\theta \]

The first example below shows the difference between the projectivity group that was already implemented in GAP and the projectivity group implemented in
\fining{}. The example also illustrates that accessing the collineation groups implemented in \fining{} is naturally associated
to the underlying incidence geometry.

\begin{example}
{\footnotesize
\begin{lstlisting}[language=GAP]
gap> g1 := PGL(5,27);
<permutation group of size 22496309500661613496614846025474560 with 2 generators>
gap> g2 := ProjectivityGroup(PG(4,27));
The FinInG projectivity group PGL(5,27)
gap> Order(g2);
22496309500661613496614846025474560
gap> g3 := CollineationGroup(PG(4,27)); 
The FinInG collineation group PGammaL(5,27)
gap> Order(g3);
67488928501984840489844538076423680
gap> GeneratorsOfGroup(g1)[1];
(20442,20467,20466,20465,20464,20463,20462,20461,20460,20459,20458,20457,
20456,20443,20455,20454,20453,20452,20451,20450,20449,20448,20447,20446,20445,
20444)(20468,21143,21116,21089,21062,21035,21008,20981,20954,20927,20900,
20873,20846,20495,20819,20792,20765,20738,20711,20684,20657,20630,20603,20576,
20549,20522)(20469,21169,21141,21113,21085,21057,21029,21001,20973,20945,
20917,20889,20861,20497,20833,20805,20777,20749,20721,20693,20665,20637,20609,
20581,20553,20525)(20470,21157,21129,21101,21073,21045,21017,20989,20961,
20933,20905,20877,20849,20496,20845,20817,20789,20761,20733,20705,20677,20649,
20621,20593,20565,20537)(20471,21144,21142,21114,21086,21058,21030,21002,
20974,20946,20918,20890,20862,20510,20821,20806,20778,20750,20722,20694,20666,
20638,20610,20582,20554,20526)(20472,21146,21117,21115,21087,21059,21031,
21003,20975,20947,20919,20891,20863,20511,20834,20794,20779,20751,20723,20695,
20667,20639,20611,20583,20555,20527)(20473,21147,21119,21090,21088,21060,
21032,21004,20976,20948,20920,20892,20864,20512,20835,20807,20767,20752,20724,
20696,20668,20640,20612,20584,20556,20528)(20474,21148,21120,21092,21063,
21061,21033,21005,20977,20949,20921,20893,20865,20513,20836,20808,20780,20740,
20725,20697,20669,20641,20613,20585,20557,20529)( [...] )
gap> GeneratorsOfGroup(g2)[1];
< a collineation: <cmat 5x5 over GF(3,3)>, F^0>
gap> GeneratorsOfGroup(g3)[1];
< a collineation: <cmat 5x5 over GF(3,3)>, F^0>
\end{lstlisting}
\caption{Projectivity groups}%
}
\end{example}

The next example shows how to construct a collineation of a projective space, and how to get the underlying matrix
and field automorphism of a collineation.

\begin{example}
{\footnotesize
\begin{lstlisting}[language=GAP]
gap> mat := [ [ Z(7^3)^115, Z(7^3)^164, Z(7^3)^189 ], 
>   [ Z(7^3)^146, Z(7^3)^279, Z(7^3)^215 ], 
>   [ Z(7^3)^58, Z(7^3)^322, Z(7^3)^260 ] ];
[ [ Z(7^3)^115, Z(7^3)^164, Z(7^3)^189 ], 
  [ Z(7^3)^146, Z(7^3)^279, Z(7^3)^215 ], 
  [ Z(7^3)^58, Z(7^3)^322, Z(7^3)^260 ] ]
gap> theta := FrobeniusAutomorphism(GF(7^3))^2;
FrobeniusAutomorphism( GF(7^3) )^2
gap> phi := ProjectiveSemilinearMap(mat,theta,GF(7^3));
< a collineation: <cmat 3x3 over GF(7,3)>, F^49>
gap> psi := Random(CollineationGroup(PG(2,2^3)));
< a collineation: <cmat 3x3 over GF(2,3)>, F^4>
gap> UnderlyingMatrix(psi);
<cmat 3x3 over GF(2,3)>
gap> Unpack(last); # to get actual matrix from cvec
[ [ Z(2^3)^6, Z(2^3), Z(2)^0 ], [ Z(2^3)^2, Z(2^3)^4, Z(2^3)^2 ], 
  [ Z(2^3)^5, Z(2^3)^3, Z(2^3)^4 ] ]
gap> FieldAutomorphism(psi);
FrobeniusAutomorphism( GF(2^3) )^2
\end{lstlisting}
\caption{Constructing collineations of a projective space}%
}
\end{example}

Note that also the projective special homography groups are implemented in GAP as permutation groups, 
these are the groups of homographies where the underlying  matrix has determinant one. 

\begin{example}
{\footnotesize
\begin{lstlisting}[language=GAP]
gap> PSL(3,49);
<permutation group of size 11072935641600 with 2 generators>
gap> SpecialHomographyGroup(PG(2,49));
The FinInG PSL group PSL(3,49)
gap> Order(last);
11072935641600
\end{lstlisting}
\caption{Special homography groups}%
}
\end{example}

Let $\calS$ be a polar space defined by the sesquilinear form $f$ on the vector space $V(d+1,q)$. 
Let $\varphi$ be a semilinear automorphism of $V(d+1,q)$. If $\varphi$ preserves the zeroes of the form $f$, that is,
 \[
 f(u^\varphi,v^\varphi) = 0 \iff f(u,v) = 0\quad \forall u,v \in V(d+1,q),
 \]
then $\varphi$ induces a collineation of $\calS$, and all collineations of $\calS$ are induced by such automorphisms of $V(d+1,q)$. 
If $\varphi$ preserves the values of the form $f$, i.e., $\forall u,v: f(u^\varphi,v^\varphi) = f(u,v)$, then the induced collineation
is called an {\em isometry}
of $\calS$. If $\varphi$ preserves the form $f$ up to a scalar, i.e., $\forall u,v: f(u^\varphi,v^\varphi) = \alpha f(u,v)$ and $\alpha$ 
is fixed, then the induced collineation is called a {\em similarity}. 
Finally, the {\em special isometries} are the isometries induced by a linear automorphism with determinant $1$, i.e., an
element of $\SL(d+1,q)$. If $\calS$ is a polar space with special
isometry group $S$, isometry group $I$, similarity group $G$ and collineation group $\Gamma$, respectively,
then $S \leqslant I \leqslant G \leqslant \Gamma$. Equalities can occur depending on the type of $\calS$. Table~\ref{groups} provides an overview. 



\begin{table}
\begin{tabular}{|c|c|c|c|c|c|}
\toprule
group/polar space & symplectic & hyperbolic & elliptic & parabolic & hermitian \\
\midrule
special isometries & $\mathrm{PSp}(d,q)$ & $\mathrm{PSO}^+(d,q)$ & $\mathrm{PSO}^-(d,q)$ & $\mathrm{PSO}(d,q)$ & $\mathrm{PSU}(d,q^2)$ \\
isometries &$\mathrm{PSp}(d,q)$ & $\mathrm{PGO}^+(d,q)$ & $\mathrm{PGO}^-(d,q)$ & $\mathrm{PGO}(d,q)$ & $\mathrm{PGU}(d,q^2)$ \\
similarities &$\mathrm{PGSp}(d,q)$ & $\mathrm{P\Delta O}^+(d,q)$ & $\mathrm{P\Delta O}^-(d,q)$ & $\mathrm{PGO}(d,q)$ & $\mathrm{PGU}(d,q^2)$ \\
collineations &$\mathrm{P\Gamma Sp}(d,q)$ & $\mathrm{P\Gamma O}^+(d,q)$ & $\mathrm{P\Gamma O}^-(d,q)$ & $\mathrm{P\Gamma O}(d,q)$ & $\mathrm{P\Gamma U}(d,q^2)$ \\
\bottomrule
\end{tabular}
\caption{Collineation groups of finite classical polar spaces}\label{groups}
\end{table}

All groups from Table~\ref{groups} are constructed in \fining{} as a collineation group of a finite classical polar space. Consistent with 
the basic philosophy of \fining{}, such a geometry is determined by a user chosen sesquilinear/quadratic form. To compute the groups, we
use the generators from \cite{KL90}, which are given with relation to a fixed sesquilinear/quadratic form, and apply the necessary 
coordinate transformations on the fly, using the package {\sf Forms}.


\begin{example}
{\footnotesize
\begin{lstlisting}[language=GAP]
gap> ps := HermitianPolarSpace(2,4);
H(2, 2^2)
gap> CollineationGroup(ps);
PGammaU(3,2^2)
gap> SpecialIsometryGroup(ps);
PSU(3,2^2)
gap> mat := IdentityMat(3,GF(9));
[ [ Z(3)^0, 0*Z(3), 0*Z(3) ], [ 0*Z(3), Z(3)^0, 0*Z(3) ], 
  [ 0*Z(3), 0*Z(3), Z(3)^0 ] ]
gap> form := QuadraticFormByMatrix(mat,GF(9));
< quadratic form >
gap> ps := PolarSpace(form);
<polar space in ProjectiveSpace(2,GF(3^2)): x_1^2+x_2^2+x_3^2=0 >
gap> g := SimilarityGroup(ps);
#I  Computing collineation group of canonical polar space...
<projective collineation group of size 720 with 2 generators>
gap> h := SimilarityGroup(ParabolicQuadric(2,9));
PGO(0,3,9)
\end{lstlisting}
\caption{Constructing collineation groups of finite classical polar spaces}%
}
\end{example}

\section{Group actions for collineation groups of Lie geometries}\label{sect:group_actions}

The general mechanism in GAP to describe the action of the elements of a group $G$ on the elements of a set $A$ is through so-called {\em action functions}. 
An action function for the group $G$ and the set $A$ is a GAP function that takes two arguments, $g \in G$ and $a \in A$, and returns the image of 
$a$ under $g$ (usually denoted as $a^g$ or $g(a)$). To explore group theoretical properties of the action of the group $G$ on the elements of a set $A$,
some generic GAP functions will compute a permutation representation of this action, which necessarily uses the action function. So for a given 
incidence geometry $\calS$, with a group of collineations $G$, \fining{} usually provides an action function for the group $G$ and the set of elements (of a given type)
of $\calS$. For Lie geometries, \fining{} provides the action function \verb|OnProjSubspaces|. The following example illustrates its use.

\begin{example}
{\footnotesize
\begin{lstlisting}[language=GAP]
gap> ps := ParabolicQuadric(4,3);
Q(4, 3)
gap> group := IsometryGroup(ps);
PGO(0,5,3)
gap> hom := ActionHomomorphism(group,Lines(ps),OnProjSubspaces);
<action homomorphism>
gap> act := Image(hom);
Group([ (1,2)(3,5)(4,7)(6,10)(9,14)(11,17)(12,18)(15,21)(16,23)(19,22)(20,24)
(25,29)(27,28)(30,34)(33,38)(36,39), (1,3,6,11)(2,4,8,13,19,14,20,25,30,35,40,
18)(5,9,15,22,27,32,37,38,34,39,29,17)(7,12)(10,16)(21,26,31,36)(23,28,33) ])
\end{lstlisting}
\caption{Action of a collineation group on projective subspaces}%
}
\end{example}

In GAP, a {\em nice monomorphism} of an (abstract) group $G$ is a permutation representation of $G$ on a suitable action domain.
Such a permutation representation makes it possible to use permutation
group algorithms to explore properties of $G$. If a nice monomorphism is available for an abstract group, a complete set of GAP functions will perform
much better than their generic counterparts for abstract groups. A typical example is the GAP command \verb|SylowSubgroup|. In the following
example the nice monomorphism of the collineation group is only computed when it becomes useful.
\begin{example}
{\footnotesize
\begin{lstlisting}[language=GAP]
gap> ps := CollineationGroup(PG(3,4));
The FinInG collineation group PGammaL(4,4)
gap> HasNiceMonomorphism(ps);
false
gap> SylowSubgroup(ps,7);
<projective collineation group of size 7>
gap> HasNiceMonomorphism(ps);
true
\end{lstlisting}
\caption{Computation of a nice monomorphism when required}%
}
\end{example}
If a collineation group $G$ is known to be a subgroup of a finite classical group, then the permutation representation of smallest
degree will be computed automatically when needed by \fining{}. 
Constructing a particular group of collineations as a subgroup of a collineation group of a polar space makes this permutation representation
available automatically for the constructed subgroup. The following example nicely illustrates this behaviour.

\begin{example}
{\footnotesize
\begin{lstlisting}[language=GAP]
gap> ps := HermitianPolarSpace(3,81);
H(3, 9^2)
gap> group := CollineationGroup(ps);
PGammaU(4,9^2)
gap> time;
2219
gap> g := Random(group);
< a collineation: <cmat 4x4 over GF(3,4)>, F^27>
gap> h := Random(group);
< a collineation: <cmat 4x4 over GF(3,4)>, F^3>
gap> group2 := Group([g,h]);
<projective collineation group with 2 generators>
gap> HasNiceMonomorphism(group2);
false
gap> Order(group2);
407194345728000
gap> time;
371559
gap> HasNiceMonomorphism(group2);
true
gap> NrMovedPoints(NiceObject(group2));
538084
gap> Size(Points(PG(3,81)));
538084
gap> group2 := Subgroup(group, [g,h]);
<projective collineation group with 2 generators>
gap> HasNiceMonomorphism(group2);
true
gap> HasNiceObject(group2);
false
gap> Order(group2);
407194345728000
gap> time;
888
gap> HasNiceObject(group2);
true
gap> NrMovedPoints(NiceObject(group2));
7300
gap> Size(Lines(ps));
7300
\end{lstlisting}
\caption{Nice monomorphisms illustrated}%
}
\end{example}

We remark that there is a byproduct of our package that would be useful to the general group theory community. The affine groups
$\agammal(d,q)$ and $\mathrm{AGL}(d,q)$, in their natural action on subspaces, are readily computable using \fining{}, and so too are correlation groups of projective spaces
(i.e., the automorphism groups of projective special linear groups). It is now a simple exercise to have the exceptional groups $G_2(q)$ and $\,^3D_4(q)$ 
acting on the subspaces in their natural linear representations.

\section{Orbits and stabilisers}

The efficient computation of orbits and stabilisers is done through the packages {\sf Orb} \cite{orb} and {\sf GenSS} \cite{genss}. 
{\sf Orb} is used for efficient orbit enumeration and the most important algorithms rely fundamentally on code for producing hash tables.
The basic idea is to enumerate an orbit by ``suborbits'' with respect to one or more subgroups, and usually requires the user
to give specific details about the putative orbit to obtain dividends in using this package over the standard method given in the GAP library.
For instance the isometry group of a polar space, acts transitively on totally isotropic subspaces of a fixed dimension, and 
we know in advance how large this orbit will be. So {\sf Orb} is ideally suited for quickly enumerating the set
of elements of a finite polar space. We also use {\sf Orb} to efficiently construct \emph{nice monomorphisms} (see Section \ref{sect:group_actions}). 
The {\sf GenSS} package stands for `Generic Schreier-Sims' and is essentially a suite of algorithms for stabiliser chains, bases, and strong
generating sets for finite groups, that are not necessarily permutation groups. Again, this removes the need for us to implement stabiliser chains
for the collineation groups in our package.

The package {\sf Orb} takes advantage of the use of the package {\sf cvec} \cite{cvec} which implements vectors and matrices over finite 
fields, replacing the existing GAP implementations. The use of {\sf cvec} is clear when asking for underlying objects of subspaces
of projective spaces, as the following example illustrates.

\begin{example}
{\footnotesize
\begin{lstlisting}[language=GAP]
gap> pg := PG(4,5);
ProjectiveSpace(4, 5)
gap> p := Random(Points(pg));
<a point in ProjectiveSpace(4, 5)>
gap> UnderlyingObject(p);
<cvec over GF(5,1) of length 5>
gap> l := Random(Lines(pg));
<a line in ProjectiveSpace(4, 5)>
gap> UnderlyingObject(l);
<cmat 2x5 over GF(5,1)>
gap> Unpack(last);
[ [ Z(5)^0, 0*Z(5), Z(5)^2, Z(5)^3, 0*Z(5) ], 
  [ 0*Z(5), Z(5)^0, Z(5)^2, Z(5)^2, Z(5)^2 ] ]
\end{lstlisting}
\caption{cvec objects as underlying objects of elements}%
}
\end{example}

Part of the generic GAP functionality is to compute orbits and stabilisers using generic algorithms. However, \fining{} implements its
own versions for Lie geometries, based on the functions in the packages {\sf Orb} and {\sf GenSS}. The following example shows
the difference in computation time between the generic GAP function \verb|Orbit|, \verb|Stabilizer| and the \fining{} function \verb|FiningOrbit|, \verb|FiningStabiliser|,
and \verb|FiningSetwiseStabiliser|. Note the rather big difference in computation time between the use of the generic function \verb|Stabilizer| (combined
with \verb|OnSets|) and \verb|FiningSetwiseStabiliser|.

\begin{example}
{\footnotesize
\begin{lstlisting}[language=GAP]
gap> ps := ParabolicQuadric(6,5);
Q(6, 5)
gap> line := Random(Lines(ps));
<a line in Q(6, 5)>
gap> group := CollineationGroup(ps);
PGammaO(7,5)
gap> orb1 := FiningOrbit(group,line,OnProjSubspaces); 
<closed orbit, 101556 points>
gap> time;
7956
gap> orb2 := Orbit(group,line,OnProjSubspaces);;
gap> time;
13075
gap> Stabilizer(group,line,OnProjSubspaces);
<projective collineation group of size 4500000000 with 5 generators>
gap> time;
152627
gap> FiningStabiliser(group,line);
<projective collineation group of size 4500000000 with 3 generators>
gap> time;
2070
gap> ps := ParabolicQuadric(4,5);
Q(4, 5)
gap> p1 := Random(Points(ps));
<a point in Q(4, 5)>
gap> p2 := Random(Points(ps));
<a point in Q(4, 5)>
gap> p3 := Random(Points(ps));
<a point in Q(4, 5)>
gap> group := CollineationGroup(ps);
PGammaO(5,5)
gap> Stabilizer(group,Set([p1,p2,p3]),OnSets);
<projective collineation group of size 72 with 4 generators>
gap> time;
113330
gap> FiningSetwiseStabiliser(group,Set([p1,p2,p3]));
#I  Computing adjusted stabilizer chain...
<projective collineation group with 4 generators>
gap> time;
18
\end{lstlisting}
\caption{Stabilizers and setwise stabilizers}%
}
\end{example}

\section{Geometry morphisms}

Let $\calS$ and $\calS'$ be two incidence geometries.  Let ${\cal E}$ be a subset of the elements of $\calS$, and let ${\cal E}'$ be a subset of the elements of $\calS'$. 
A map $\varphi: {\cal E} \rightarrow {\cal E}'$ is called a geometry morphism from $\calS$ to $\calS'$ if it preserves incidence:
\[ \forall e,f \in {\cal E}: \varphi(e) \I' \varphi(f) \iff e \I f \,.\]
Note that $\cE$, $\cE'$ respectively, might be any subset of the elements of $\calS$, $\calS'$ respectively. Let $G$, respectively $G'$ be 
the automorphism group of $\calS$, respectively $\calS'$. A map $\theta: G \rightarrow G'$ is an {\em intertwiner} for the geometry morphism
$\varphi$ if
\[\forall e \in \cE, \forall g \in G: \varphi(e^g) = \varphi(e)^{\theta(g)}\,.\]
The definition of a geometry morphism is very general. A geometry morphism might have the complete set of elements of an incidence geometry
as domain and might be type preserving.  Geometry morphisms with particular properties occur frequently in incidence geometry. 
Typical geometry morphisms are embeddings, dualities, the Klein correspondence, morphisms induced by field reduction, etc.  All geometry
morphisms currently available in \fining{} are geometry morphisms of Lie geometries. We give some particular examples. 

We've seen that classical polar spaces can be constructed through a user-defined form. The isomorphism between two polar spaces
defined from a similar but different form is present as a geometry morphism. Clearly this morphism has the complete set of elements
of the polar spaces as domain/co-domain, and is type preserving.
\begin{example}
{\footnotesize
\begin{lstlisting}[language=GAP]
gap> mat1 := [[1,0,0],[0,0,1],[0,0,0]]*Z(5)^0;
[ [ Z(5)^0, 0*Z(5), 0*Z(5) ], [ 0*Z(5), 0*Z(5), Z(5)^0 ], 
  [ 0*Z(5), 0*Z(5), 0*Z(5) ] ]
gap> form1 := QuadraticFormByMatrix(mat1,GF(5));
< quadratic form >
gap> conic1 := PolarSpace(form1);
gap> mat2 := [[0,0,-1],[0,1,0],[0,0,0]]*Z(5)^0; 
[ [ 0*Z(5), 0*Z(5), Z(5)^2 ], [ 0*Z(5), Z(5)^0, 0*Z(5) ], 
  [ 0*Z(5), 0*Z(5), 0*Z(5) ] ]
gap> form2 := QuadraticFormByMatrix(mat2,GF(5));
< quadratic form >
gap> conic2 := PolarSpace(form2);
<polar space in ProjectiveSpace(2,GF(5)): -x_1*x_3+x_2^2=0 >
gap> iso := IsomorphismPolarSpaces(conic1,conic2);
#I  No intertwiner computed. One of the polar spaces must 
have a collineation group computed
<geometry morphism from <Elements of Q(2, 5): x_1^2+x_2*x_3=0> 
to <Elements of Q(2, 5): -x_1*x_3+x_2^2=0>>
\end{lstlisting}
\caption{Isomorphic polar spaces}%
}
\end{example}

Example~\ref{ex:subspace_embedding}, dealing with subpace embeddings, is also among the most obvious ones. It is demonstrated how to embed a projective plane
in a projective space as a user chosen plane of the projective space, and similarly, how to embed a conic in an elliptic quadric as a user chosen 
plane section.

\begin{example}
{\footnotesize
\begin{lstlisting}[language=GAP]
gap> mat := [ [ Z(5)^0, 0*Z(5), Z(5), 0*Z(5), 0*Z(5) ], 
>   [ 0*Z(5), Z(5)^0, Z(5)^3, 0*Z(5), 0*Z(5) ], 
>   [ 0*Z(5), 0*Z(5), 0*Z(5), Z(5)^0, Z(5)^2 ] ];;
gap> plane := VectorSpaceToElement(PG(4,5),mat);
<a plane in ProjectiveSpace(4, 5)>
gap> em := NaturalEmbeddingBySubspace(PG(2,5),PG(4,5),plane);
<geometry morphism from <All elements of ProjectiveSpace(2, 
5)> to <All elements of ProjectiveSpace(4, 5)>>
gap> ps := EllipticQuadric(5,11);
Q-(5, 11)
gap> mat := [ [ Z(11)^0, 0*Z(11), 0*Z(11), 0*Z(11), 0*Z(11), 0*Z(11) ], 
>   [ 0*Z(11), Z(11)^0, 0*Z(11), 0*Z(11), 0*Z(11), 0*Z(11) ], 
>   [ 0*Z(11), 0*Z(11), Z(11)^0, 0*Z(11), Z(11)^0, Z(11)^0 ] ];;
gap> plane := VectorSpaceToElement(PG(5,11),mat);
<a plane in ProjectiveSpace(5, 11)>
gap> TypeOfSubspace(ps,plane);
"parabolic"
gap> conic := ParabolicQuadric(2,11);
Q(2, 11)
gap> em := NaturalEmbeddingBySubspace(conic,ps,plane);
<geometry morphism from <Elements of Q(2, 11)> to <Elements of Q-(5, 11)>>
\end{lstlisting}
\caption{Subspace embeddings}\label{ex:subspace_embedding}%
}
\end{example}

The Klein correspondence is a geometry morphism that has the set of lines of $\pg(3,q)$ as domain and the set of points of $\q^+(5,q)$ as
co-domain. This morphism is clearly not type preserving. Example~\ref{ex:klein_dualities} also shows some of the derived dualites 
of the particular finite generalised quadrangles.

\begin{example}
{\footnotesize
\begin{lstlisting}[language=GAP]
gap> map := KleinCorrespondence(7);
<geometry morphism from <lines of ProjectiveSpace(3, 7)> 
to <points of Q+(5, 7): x_1*x_6+x_2*x_5+x_3*x_4=0>>
gap> p := Random(Lines(PG(3,7)));
<a line in ProjectiveSpace(3, 7)>
gap> p^map;
<a point in Q+(5, 7): x_1*x_6+x_2*x_5+x_3*x_4=0>
gap> duality := NaturalDuality(ParabolicQuadric(4,3),SymplecticSpace(3,3));
<geometry morphism from <Elements of Q(4, 3)> to <Elements of W(3, 3)>>
gap> w := AmbientGeometry(Range(duality));
W(3, 3)
gap> l := Random(Points(w));
<a point in W(3, 3)>
gap> PreImageElm(duality,l);
<a line in Q(4, 3)>
\end{lstlisting}
\caption{The Klein correspondence}\label{ex:klein_dualities}%
}
\end{example}

Field extensions are the basic ingredient for some important embeddings. Typically the associated geometry morphism has the complete set of
elements of a geometry as domain, and a subset of the elements of a geometry as co-domain. We give some examples: the embedding of 
$\pg(n,q)$ in $\pg(n,q^t)$ (illustrated using $\pg(2,3)$ embedded in $\pg(2,27)$), and the embedding of $\q^-(5,q)$ in $\q^+(5,q^2)$ 
(illustrated using $\q^-(5,2)$ embedded in $\q^+(5,4)$). Note
that the second example requires the use of coordinate transformations, since considering the form of the standard elliptic quadric $\q^-(5,q)$ over $\gf(q^2)$
is not necessarily the form used for the standard hyperbolic quadric $\q^+(5,q^2)$, but this calculation is done on the fly. Afterwards, we 
compute the Klein correspondence from $\pg(3,4)$ to $\q^+(5,4)$, compute the lines of $\pg(3,4)$ corresponding with the points of $\q^-(5,2)$ embedded
in $\q^+(5,4)$, compute the set of points on at least one of the lines of this line set, and then the intersection numbers of the lines of $\pg(3,4)$ with
this point set. The intersection numbers are sufficient to conclude that the point set is the set of points of a Hermitian polar space $\h(3,4)$, (by
a result of \cite{DWS2010}), and illustrate nicely how the duality between the $\q^-(5,q)$ and $\h(3,q^2)$ is obtained through the Klein correspondence.

\begin{example}
{\footnotesize
\begin{lstlisting}[language=GAP]
gap> em := NaturalEmbeddingBySubfield(PG(2,3),PG(2,3^3));
<geometry morphism from <All elements of ProjectiveSpace(2, 3)> 
to <All elements of ProjectiveSpace(2, 27)>>
gap> Collected(List(Lines(PG(2,3^3)),x->Number(Points(PG(2,3)),y->y^em in x)));
[ [ 0, 432 ], [ 1, 312 ], [ 4, 13 ] ]
gap> em := NaturalEmbeddingBySubfield(PG(2,3),PG(2,3^2));
<geometry morphism from <All elements of ProjectiveSpace(2, 3)> 
to <All elements of ProjectiveSpace(2, 9)>>
gap> Collected(List(Lines(PG(2,3^2)),x->Number(Points(PG(2,3)),y->y^em in x)));
[ [ 1, 78 ], [ 4, 13 ] ]
gap> e := EllipticQuadric(5,2);
Q-(5, 2)
gap> h := HyperbolicQuadric(5,4);
Q+(5, 4)
gap> em := NaturalEmbeddingBySubfield(e,h);
#I  No intertwiner computed. <geom1> must have a collineation group computed
<geometry morphism from <Elements of Q-(5, 2)> to <Elements of Q+(5, 4)>>
gap> klein := KleinCorrespondence(h);
<geometry morphism from <lines of ProjectiveSpace(3, 4)> to <points of Q+(5, 
4)>>
gap> lines := List(Points(e),x->PreImageElm(klein,x^em));;
gap> pointset := Union(List(lines,x->List(Points(x))));;
gap> Collected(List(Lines(PG(3,4)),x->Number(pointset,y->y in x)));
[ [ 1, 90 ], [ 3, 240 ], [ 5, 27 ] ]
gap> Collected(List(Planes(PG(3,4)),x->Number(pointset,y->y in x)));
[ [ 9, 40 ], [ 13, 45 ] ]
\end{lstlisting}
\caption{Subfield embeddings}%
}
\end{example}

Field reduction is a way to embed low rank Lie geometries over a large field into high rank Lie geometries over a subfield. An obvious example 
is the embedding of, for example, $\pg(2,q^2)$ into $\pg(5,q)$. The corresponding geometry morphism will have the points and lines of $\pg(2,q^2)$ as domain, and
a subset of the lines and $3$-dimensional subspaces of $\pg(5,q)$ as co-domain. Similarly, embeddings of polar spaces are possible. Let $\calS$ be a classical 
polar space with associated form $f$ acting on the vector space $V(r,q^t)$. Let $\mathrm{Tr}$ be the trace map from $\gf(q^t)$ to $\gf(q)$. Note that choosing
a basis for $\gf(q^t)$ as a $t$-dimensional vector space over $\gf(q)$, induces a bijection $\Phi$ between the vectors of $V(r,q^t)$ and the vectors of $V(rt,q)$. 
Define, for a chosen $\alpha \in \gf(q^t) \setminus \{0\}$, the map $\mathrm{Tr}_{\alpha}: \gf(q^t) \rightarrow \gf(q)$: $x \mapsto \mathrm{Tr}(\alpha x)$.
Then clearly $\mathrm{Tr}_{\alpha} \circ f \circ \Phi^{-1}$ is a form on $V(rt,q)$, defining a classical polar space $\calS'$. An element of $\calS$ is embedded
in $\calS'$ by blowing up the underlying vector space. Depending on the choice of $\alpha$, some conditions on $q$ and the type of the form $f$, several embeddings
are possible. A complete overview (including all proofs), is found in \cite{LV13}. Note that these embeddings are very useful in constructing particular substructures
of polar spaces over small fields, as the examples below shows by constructing a so-called hermitian $1$-system of $\w(5,2)$.

\begin{example}
{\footnotesize
\begin{lstlisting}[language=GAP]
gap> h := HermitianPolarSpace(2,4);
H(2, 2^2)
gap> w := SymplecticSpace(5,2);
W(5, 2)
gap> em := NaturalEmbeddingByFieldReduction(h,w);
#I  These polar spaces are suitable for field reduction
<geometry morphism from <Elements of H(2, 2^2)> to <Elements of W(5, 2)>>
gap> List(Points(h),x->x^em);
[ <a line in W(5, 2)>, <a line in W(5, 2)>, <a line in W(5, 2)>, 
  <a line in W(5, 2)>, <a line in W(5, 2)>, <a line in W(5, 2)>, 
  <a line in W(5, 2)>, <a line in W(5, 2)>, <a line in W(5, 2)> ]
\end{lstlisting}
\caption{Embeddings by field reduction}%
}
\end{example}

\section{Enumerators for elements of classical polar spaces}

Let $\cC$ be the collection (in the GAP sense) of elements of a classical polar space $\calS$. Note that a collection in the GAP sense 
is an abstract object, which just {\em represents} all elements (of a given type) rather than contains them. At creation, there is no computation needed. An enumerator
for a collection is a GAP-object that allows one to compute the $i$th element of the collection $\cC$ without computing any other element 
of $\cC$ and vice versa, i.e., given an element $e \in \cC$, the enumerator is able to compute its position in the list $\cC$ without enumerating
any of the elements of $\cC$. Given a collection $\cC$ for which the number of elements is known, an enumerator enables us to 
efficiently compute random elements of $\cC$ if no other information is known. The algorithms used for this part of \fining{} were developed 
specifically for this purpose.
\begin{example}
{\footnotesize
\begin{lstlisting}[language=GAP]
gap> lines := Lines( ParabolicQuadric(6,3) );
<lines of Q(6, 3)>
gap> enum := Enumerator( lines );
EnumeratorOfSubspacesOfClassicalPolarSpace( <lines of Q(6, 3)> )
gap> s := Size(enum);
3640
gap> n := Random([1..s]);
3081
gap> l := enum[n];
<a line in Q(6, 3)>
gap> Position(enum, l);
3081
\end{lstlisting}
\caption{Enumerator of the points of a polar space}%
}
\end{example}

\begin{example}
{\footnotesize
\begin{lstlisting}[language=GAP]
gap> ps := ParabolicQuadric(16,7^4);
Q(16, 2401)
gap> pts := Points(ps);
<points of Q(16, 2401)>
gap> Size(pts);
508233536514931541724405776067904925314839705888016
gap> Random(pts);
<a point in Q(16, 2401)>
gap> time;
565
\end{lstlisting}
\caption{Enumerators and their use in \fining{}}%
}
\end{example}

\section{Coset geometries}
Consider an incidence geometry $\mathcal{S}$ together with a group $G$
of automorphisms which is transitive on the set of maximal flags. For
such a maximal flag $C=\{c_1,c_2,\ldots ,c_n\}$ we write $G_i$ for the
stabiliser of the element $c_i$ in $G$ where $1\leq i\leq n$. We will
assume that the index $i$ corresponds to the type of $c_i$ in $C$. It
is then straightforward to see that the elements of type $i$ in
$\mathcal{S}$ can be identified with the (left) cosets of $G_i$ in
$G$. In this representation, the incidence relation of
$\mathcal{S}$ is translated to nonempty intersection of cosets,
i.e., two elements identified with $gG_i$ and $hG_j$ respectively are
incident if and only if $gG_i\cap hG_j\neq \varnothing$. The action of
$G$ on $\mathcal{S}$ is carried to left multiplication: an element
identified with $hG_i$ is mapped by $g\in G$ to $ghG_i$. 

We can now reverse this construction and start with a group $G$
together with a set $\{G_1,G_2,\ldots ,G_n\}$ of subgroups of
$G$. From this we define an incidence structure whose elements are
simply the cosets of the given subgroups and incidence is given by
nonempty intersection of cosets. The group $G$ acts by left
multiplication. In general this structure is not an incidence geometry
but one can provide necessary conditions for this (see~\cite{BC13}). 

Once an incidence geometry is represented as a coset geometry, GAP's
built-in machinery for (permutation) groups can be used to analyse
this geometry. 

\begin{example}
{\footnotesize
\begin{lstlisting}[language=GAP]
gap> g:=SymmetricGroup(8);;
gap> parabolics:=[];;
gap> parabolics[1]:=Stabilizer(g,1);; parabolics[2]:=Stabilizer(g,2);;
gap> parabolics[3]:=Stabilizer(g,3);;
gap> parabolics[4]:=Stabilizer(g,[1,2,3,4],OnSets);;
gap> parabolics[5]:=Stabilizer(g,[1,2,3,4,5],OnSets);;
gap> parabolics[6]:=Stabilizer(g,6);; parabolics[7]:=Stabilizer(g,7);;
gap> cg:=CosetGeometry(g,parabolics);
CosetGeometry( SymmetricGroup( [ 1 .. 8 ] ) )
gap> IsFirmGeometry(cg);
true
gap> IsThinGeometry(cg);
true
gap> IsThickGeometry(cg);
false
gap> truncation:=CosetGeometry(g,parabolics{[1..5]});
CosetGeometry( SymmetricGroup( [ 1 .. 8 ] ) )
gap> IsFirmGeometry(truncation);
true
gap> IsThinGeometry(truncation);
false
gap> IsThickGeometry(truncation);
false
gap> truncation2:=CosetGeometry(g,parabolics{[4,5]});
CosetGeometry( SymmetricGroup( [ 1 .. 8 ] ) )
gap> IsFirmGeometry(truncation2);
true
gap> IsThinGeometry(truncation2);
false
gap> IsThickGeometry(truncation2);
true
\end{lstlisting}
\caption{A coset geometry}\label{ex:cosetgeometry}%
}
\end{example}

One of the concepts summarizing the properties of an
incidence geometry (on which a group acts flag-transitively) is that
of a \emph{diagram}. We provide functionality to compute the diagram
of a coset geometry and also to actually draw the diagram. Since a
diagram is in fact a labeled graph, drawing is done by producing
some code which can be processed by the open source package
\verb+Graphviz+ \cite{Graphviz}. If the
package (which is not a GAP package) is installed on your computer, it
will be invoked automatically and the result will be a PostScript file
with a picture of the diagram. 

\begin{example}
{\footnotesize
\begin{lstlisting}[language=GAP]
gap> diag:=DiagramOfGeometry(truncation);
gap> DrawDiagram(diag, "truncimg");
gap> Exec("gv truncimg.ps");
\end{lstlisting}
\caption{Creating the diagram of a coset geometry (see also Example~\ref{ex:cosetgeometry})}%
}
\end{example}

\begin{figure}
\begin{center}
\includegraphics[width=.5\textwidth]{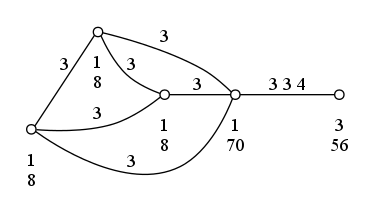}
\end{center}
\caption{The diagram of the coset geometry from Example~\ref{ex:cosetgeometry}}
\end{figure}
Interaction is also provided with the \textsf{GRAPE} package which is
very useful in studying a geometry as a (coloured) graph. Using the
\verb+IncidenceGraph+ command one obtains the incidence graph of an
incidence structure $\mathcal{S}$. This can then be used, for example,
to obtain the full automorphism group of the structure, using the very
fast \verb+nauty+ package included with \textsf{GRAPE}. Some shortcut
commands have been included for that.

\begin{example}
{\footnotesize
\begin{lstlisting}[language=GAP]
gap> g := PSL(2,11);;
gap> g1 := Group([ (1,2,3)(4,8,12)(5,10,9)(6,11,7), 
(1,2)(3,4)(5,12)(6,11)(7,10)(8,9) ]);;
gap> g2 := Group([ (1,2,7)(3,9,4)(5,11,10)(6,8,12), 
(1,2)(3,4)(5,12)(6,11)(7,10)(8,9) ]);;
gap> g3 := Group([ (1,2,11)(3,8,7)(4,9,5)(6,10,12), 
(1,2)(3,12)(4,11)(5,10)(6,9)(7,8) ]);;
gap> g4 := Group([ (1,2,11)(3,8,7)(4,9,5)(6,10,12), 
(1,2)(3,10)(4,9)(5,8)(6,7)(11,12) ]);;
gap> cg:=CosetGeometry(g,[g1,g2,g3,g4]);
CosetGeometry( Group( [ ( 3,11, 9, 7, 5)( 4,12,10, 8, 6), 
 ( 1, 2, 8)( 3, 7, 9)( 4,10, 5)( 6,12,11) ] ) )
gap> aut:=AutGroupIncidenceStructureWithNauty(cg);
<permutation group with 4 generators>
gap> StructureDescription(aut);
"PSL(2,11)"
gap> cor:=CorGroupIncidenceStructureWithNauty(cg);
<permutation group with 5 generators>
gap> StructureDescription(cor);
"C2 x PSL(2,11)"
\end{lstlisting}
\caption{Incidence geometries and graphs}%
}
\end{example}
\section{Generalised polygons}

A generalised $n$-gon is a point/line geometry whose incidence graph is bipartite of diameter $n$ and girth $2n$. Generalised 
polygons were introduced by J. Tits in \cite{Tits59}. Projective planes are examples of generalised $3$-gons and classical 
polar spaces of rank 2 are examples of generalised $4$-gons, also called generalised quadrangles. Non-classical examples 
of projective planes and generalised quadrangles are well-known. By the famous theorem of 
Feit and Higman, a generalised $n$-gon which has at least three points on every line, must have $n \in \{2,3,4,6,8\}$. 
Two classes of examples of generalised hexagons are the so-called split Cayley hexagons and the twisted triality hexagons, 
which are both derived from a triality of the Lie geometries of type $D_4$, which are in the finite case subgeometries of the hyperbolic quadric $\q^+(7,q)$. 
Generalised octagons seem to be rare. In the finite case, only the so-called Ree-Tits octagon (together with its dual) is known.
Generalised polygons occur in several ways and one particular generalised polygon can often be represented in different ways. 
The framework for generalised polygons in \fining{} is generic. It provides the user with a set of operations that are always applicable
to any generalised polygon (and its elements), independent of its representation. The classical examples such as classical projective planes, 
classical generalised quadrangles and the classical generalised hexagons are represented as Lie geometries, but are also
generalised polygons in \fining{} and as such all generic operations for generalised polygons can be used. This generic framework also
includes operations that deal with the incidence graph of the generalised polygon.
%
%

We give some examples of the different possibilities to construct generalised polygons in \fining{}. In Example~\ref{ex:gp_blocks} a projective 
plane of order $4$ is constructed by representing every line as a set of points. The points may be any set of objects, 
as long as the set of blocks is a set of subsets of the point set. In Example~\ref{ex:tits} the construction of the generalised quadrangle 
$T_2({\mathcal{O}})$ by using elements of a suitable projective space, is demonstrated. Finally, Example~\ref{ex:splitcayley} 
demonstrates the construction of the Split Cayley hexagon. Note that the construction operations used in the first two examples
\verb|GeneralisedPolygonByBlocks| and \verb|GeneralisedPolygonByElements| check whether the input really yields a generalised polygon. 
\begin{example}
{\footnotesize
\begin{lstlisting}[language=GAP]
gap> blocks := [ 
>   [ 1, 2, 3, 4, 5 ], [ 1, 6, 7, 8, 9 ], [ 1, 10, 11, 12, 13 ],
>   [ 1, 14, 15, 16, 17 ], [ 1, 18, 19, 20, 21 ], [ 2, 6, 10, 14, 18 ], 
>   [ 2, 7, 11, 15, 19 ], [ 2, 8, 12, 16, 20 ], [ 2, 9, 13, 17, 21 ], 
>   [ 3, 6, 11, 16, 21 ], [ 3, 7, 10, 17, 20 ], [ 3, 8, 13, 14, 19 ], 
>   [ 3, 9, 12, 15, 18 ], [ 4, 6, 12, 17, 19 ], [ 4, 7, 13, 16, 18 ], 
>   [ 4, 8, 10, 15, 21 ], [ 4, 9, 11, 14, 20 ], [ 5, 6, 13, 15, 20 ], 
>   [ 5, 7, 12, 14, 21 ], [ 5, 8, 11, 17, 18 ], [ 5, 9, 10, 16, 19 ] ];;
gap> gp := GeneralisedPolygonByBlocks( blocks );
<projective plane order 4>
\end{lstlisting}
\caption{A projective plane constructed by its blocks}\label{ex:gp_blocks}%
}
\end{example}
\begin{example}
{\footnotesize
\begin{lstlisting}[language=GAP]
gap> conic := Set(Points(ParabolicQuadric(2,4)));
[ <a point in Q(2, 4)>, <a point in Q(2, 4)>, <a point in Q(2, 4)>, 
  <a point in Q(2, 4)>, <a point in Q(2, 4)> ]
gap> pg := PG(3,4);
ProjectiveSpace(3, 4)
gap> hyp :=  HyperplaneByDualCoordinates(pg,[1,0,0,0]*Z(q)^0);
<a plane in ProjectiveSpace(3, 4)>
gap> em := NaturalEmbeddingBySubspace(PG(2,4),pg,hyp);
<geometry morphism from <All elements of ProjectiveSpace(2, 4)> 
to <All elements of ProjectiveSpace(3, 4)>>
gap> conic_at_infinity := List(conic,x->x^em);;
gap> group := CollineationGroup(pg);
The FinInG collineation group PGammaL(4,4)
gap> stab := FiningSetwiseStabiliser(group,conic_at_infinity);
#I  Computing adjusted stabilizer chain...
<projective collineation group with 6 generators>
gap> points1 := Set(Filtered(Points(pg),x->not x in hyp));;
gap> tangents := List(conic,x->TangentSpace(x)^em);
[ <a line in ProjectiveSpace(3, 4)>, <a line in ProjectiveSpace(3, 4)>, 
  <a line in ProjectiveSpace(3, 4)>, <a line in ProjectiveSpace(3, 4)>, 
  <a line in ProjectiveSpace(3, 4)> ]
gap> planes := List(tangents,x->Filtered(Planes(x),y->not y in hyp));;
gap> points2 := Union(planes);;
gap> points3 := [hyp];
[ <a plane in ProjectiveSpace(3, 4)> ]
gap> linesa := Union(List(conic_at_infinity,x->Filtered(Lines(x),y->not y in hyp)));;
gap> linesb := Set(conic_at_infinity);;
gap> pts := Union(points1,points2,points3);;
gap> lns := Union(linesa,linesb);;
gap> inc := \*;
<Operation "*">
gap> gp := GeneralisedPolygonByElements(pts,lns,inc,stab,\^);
<generalised quadrangle of order [ 4, 4 ]>
\end{lstlisting}
\caption{A generalised quadrangle $T_2(O)$}\label{ex:tits}
}
\end{example}
\begin{example}
{\footnotesize
\begin{lstlisting}[language=GAP]
gap> gp := SplitCayleyHexagon(3);
H(3)
\end{lstlisting}
\caption{A classical generalised hexagon}\label{ex:splitcayley}%
}
\end{example}

Example~\ref{ex:ree_tits_elemenst} shows some particular operations on (the elements of) generalised polygons. 
We use the Ree-Tits octagon as example. Note that the construction of this octagon is done by 
\verb|GeneralisedPolygonByElements| based on the construction as a coset geometry, also demonstrating
that the construction operations for generalised polygons are generic.
\begin{example}
{\footnotesize
\begin{lstlisting}[language=GAP]
gap> LoadPackage( "AtlasRep" );
true
gap> titsgroup:=AtlasGroup("2F4(2)'");
<permutation group of size 17971200 with 2 generators>
gap> g1:=AtlasSubgroup(titsgroup,3);
<permutation group of size 10240 with 2 generators>
gap> g2:=AtlasSubgroup(titsgroup,5);
<permutation group of size 6144 with 2 generators>
gap> conj:=ConjugacyClassSubgroups(titsgroup,g1);;
gap> # Now look for the conjugate of g1 with maximal intersection
gap> g1:=First(conj, sg -> Size(Intersection(sg,g2))=2048);
<permutation group of size 10240 with 2 generators>
gap> cg:=CosetGeometry(titsgroup,[g1,g2]);;
gap> pts := Set(ElementsOfIncidenceStructure(cg,1));;
gap> lines := Set(ElementsOfIncidenceStructure(cg,2));;
gap> gp := GeneralisedPolygonByElements(pts,lines,\*,titsgroup,OnCosetGeometryElement);
<generalised octagon of order [ 2, 4 ]>
gap> Order(gp);
[ 2, 4 ]
gap> p := Random(Points(gp));
<a point in <generalised octagon of order [ 2, 4 ]>>
gap> lines := Lines(p);
<shadow lines in <generalised octagon of order [ 2, 4 ]>>
gap> l := Random(lines);
<a line in <generalised octagon of order [ 2, 4 ]>>
gap> DistanceBetweenElements(p,l);
1
gap> m := Random(Lines(gp));
<a line in <generalised octagon of order [ 2, 4 ]>>
gap> DistanceBetweenElements(m,l);
8
\end{lstlisting}
\caption{The Ree-Tits octagon from a coset geometry}\label{ex:ree_tits_elemenst}%
}
\end{example}
The third example shows these particular examples on the classical generalised polygons.

\begin{example}[h!]
{\footnotesize
\begin{lstlisting}[language=GAP]
gap> h := TwistedTrialityHexagon(2^3);
T(8, 2)
gap> p := Random(Points(h));
#I  Computing nice monomorphism...
#I  Found permutation domain...
<a point in T(8, 2)>
gap> lines := Lines(p);
<shadow lines in T(8, 2)>
gap> l := Random(lines);
<a line in T(8, 2)>
gap> DistanceBetweenElements(p,l);
1
gap> m := Random(Lines(h));
<a line in T(8, 2)>
gap> DistanceBetweenElements(m,l);
6
\end{lstlisting}
\caption{Elements of the twisted triality hexagon}%
}
\end{example}

\section{Algebraic varieties}

In \fining{} some basic functionality for algebraic varieties defined over finite fields is implemented. 
The algebraic varieties in \fining{} are defined by a list of multivariate polynomials over a finite field, and an ambient geometry. 
This ambient geometry is either a projective space, in which case the algebraic variety is called a projective variety, or an affine geometry, 
in which case the algebraic variety is called an affine variety. 
There is an iterator for the set of points that lie on the variety and the user is enabled to randomly choose one of these points.

\begin{example}
{\footnotesize
\begin{lstlisting}[language=GAP]
gap> r:=PolynomialRing(GF(5),6);
GF(5)[x_1,x_2,x_3,x_4,x_5,x_6]
gap> f1:=Sum(List(IndeterminatesOfPolynomialRing(r),x->x^2));
x_1^2+x_2^2+x_3^2+x_4^2+x_5^2+x_6^2
gap> f2:=r.1*r.2+r.3*r.4+r.5*r.6;
x_1*x_2+x_3*x_4+x_5*x_6
gap> var:=AlgebraicVariety(PG(5,5),[f1,f2]);
Projective Variety in ProjectiveSpace(5, 5)
gap> Points(var);        
<points of Projective Variety in ProjectiveSpace(5, 5)>
gap> Size(Points(var));                     
156
gap> Random(Points(var));                                    
<a point in ProjectiveSpace(5, 5)>
\end{lstlisting}
\caption{Constructing an algebraic variety}%
}
\end{example}

Specific methods are available for Hermitian varieties and quadrics, as well as a method \verb|PolarSpace|
which returns the associated polar space consistent with the other functionality for polar spaces defined from forms.

\begin{example}
{\footnotesize
\begin{lstlisting}[language=GAP]
gap> hv:=HermitianVariety(3,9);                          
Hermitian Variety in ProjectiveSpace(3, 9)
gap> hp:=PolarSpace(hv);
<polar space in ProjectiveSpace(3,GF(3^2)): x_1^4+x_2^4+x_3^4+x_4^4=0 >
gap> DefiningListOfPolynomials(hv);
[ x_1^4+x_2^4+x_3^4+x_4^4 ]
gap> hv:=HermitianVariety(3,9);              
Hermitian Variety in ProjectiveSpace(3, 9)
gap> DefiningListOfPolynomials(hv);
[ x_1^4+x_2^4+x_3^4+x_4^4 ]
gap> hp:=PolarSpace(hv);
<polar space in ProjectiveSpace(3,GF(3^2)): x_1^4+x_2^4+x_3^4+x_4^4=0 >
gap> PolarSpaceType(hp);
"hermitian"
gap> IsGeneralisedQuadrangle(hp);  
true
gap> Order(hp);
[ 9, 3 ]
\end{lstlisting}
\caption{A Hermitian variety}%
}
\end{example}

The package \fining{} also contains the Veronese varieties, the Segre varieties and the Grassmann varieties; 
three classical projective varieties. These varieties have an associated geometry morphism ({\tt VeroneseMap}, 
{\tt SegreMap}, and {\tt GrassmannMap}) and  \fining{} also provides some general functionality for these. Example~\ref{ex:veronese}
shows this for the Veronese variety in a five-dimensional projective space, by constructing the set of 
planes intersecting the variety in a conic.

\begin{example}
{\footnotesize
\begin{lstlisting}[language=GAP]
gap> vv:=VeroneseVariety(PG(2,3)); 
Veronese Variety in ProjectiveSpace(5, 3)
gap> DefiningListOfPolynomials(vv);
[ -x_1*x_4+x_2^2, -x_1*x_6+x_3^2, -x_4*x_6+x_5^2, x_1*x_5-x_2*x_3 ]
gap> vm:=VeroneseMap(vv);
Veronese Map of <points of ProjectiveSpace(2, 3)>
gap> Range(vm);
<points of Veronese Variety in ProjectiveSpace(5, 3)>
gap> vvpts:=Set(Points(vv));;
gap> linesets:=List(Lines(PG(2,3)),l->List(Points(l)));;
gap> image_list:=List(linesets,list->Span(ImagesSet(vm,list)));
[ <a plane in ProjectiveSpace(5, 3)>, <a plane in ProjectiveSpace(5, 3)>, 
  <a plane in ProjectiveSpace(5, 3)>, <a plane in ProjectiveSpace(5, 3)>, 
  <a plane in ProjectiveSpace(5, 3)>, <a plane in ProjectiveSpace(5, 3)>, 
  <a plane in ProjectiveSpace(5, 3)>, <a plane in ProjectiveSpace(5, 3)>, 
  <a plane in ProjectiveSpace(5, 3)>, <a plane in ProjectiveSpace(5, 3)>, 
  <a plane in ProjectiveSpace(5, 3)>, <a plane in ProjectiveSpace(5, 3)>, 
  <a plane in ProjectiveSpace(5, 3)> ]
gap> List(image_list,plane->Number(Points(plane),x->x in vvpts));
[ 4, 4, 4, 4, 4, 4, 4, 4, 4, 4, 4, 4, 4 ]
gap> vvstab:=FiningSetwiseStabiliser(HomographyGroup(PG(5,3)),vvpts);
#I  Computing adjusted stabilizer chain...
<projective collineation group with 9 generators>
gap> plane_orbs:=FiningOrbits(vvstab,Planes(PG(5,3)));
1%..3%..4%..5%..6%..7%..12%..14%..23%..28%..29%..38%..39%..
43%..52%..60%..61%..63%..66%..68%..77%..82%..90%..99%..
100%..[ <closed orbit, 13 points>, <closed orbit, 13 points>, 
  <closed orbit, 104 points>, <closed orbit, 468 points>, 
  <closed orbit, 468 points>, <closed orbit, 468 points>, 
  <closed orbit, 468 points>, <closed orbit, 234 points>, 
  <closed orbit, 468 points>, <closed orbit, 1404 points>, 
  <closed orbit, 936 points>, <closed orbit, 2808 points>, 
  <closed orbit, 1872 points>, <closed orbit, 351 points>, 
  <closed orbit, 2808 points>, <closed orbit, 234 points>, 
  <closed orbit, 351 points>, <closed orbit, 1404 points>, 
  <closed orbit, 2808 points>, <closed orbit, 2808 points>, 
  <closed orbit, 468 points>, <closed orbit, 702 points>, 
  <closed orbit, 936 points>, <closed orbit, 702 points>, 
  <closed orbit, 2808 points>, <closed orbit, 1872 points>, 
  <closed orbit, 2808 points>, <closed orbit, 2808 points>, 
  <closed orbit, 144 points>, <closed orbit, 144 points> ]
gap> Filtered(plane_orbs,o->Number(Points(o[1]),x->x in vvpts)=4);
[ <closed orbit, 13 points> ]
gap> conicplanes:=last[1];
<closed orbit, 13 points>
gap> ForAll(image_list,plane->plane in conicplanes);
true
\end{lstlisting}
\caption{A Veronese variety}\label{ex:veronese}%
}
\end{example}


\begin{thebibliography}{10}

\bibitem{fining}
{\sc J.~Bamberg, A.~Betten, P.~Cara, J.~De~Beule, M.~Lavrauw, and
  M.~Neunh\"offer}, {\em {FinInG -- Finite Incidence Geometry, Version 1.3.3}},
  2016.

\bibitem{forms}
{\sc J.~Bamberg and J.~De~Beule}, {\em {FORMS -- Sesquilinear and Quadratic,
  Version 1.2.2}}, 2011.

\bibitem{magma}
{\sc W.~Bosma, J.~Cannon, and C.~Playoust}, {\em The {M}agma algebra system.
  {I}. {T}he user language}, J. Symbolic Comput., 24 (1997), pp.~235--265.
\newblock Computational algebra and number theory (London, 1993).

\bibitem{BC13}
{\sc F.~Buekenhout and A.~M. Cohen}, {\em Diagram geometry}, vol.~57 of
  Ergebnisse der Mathematik und ihrer Grenzgebiete. 3. Folge. A Series of
  Modern Surveys in Mathematics [Results in Mathematics and Related Areas. 3rd
  Series. A Series of Modern Surveys in Mathematics], Springer, Heidelberg,
  2013.
\newblock Related to classical groups and buildings.

\bibitem{DWS2010}
{\sc S.~De~Winter and J.~Schillewaert}, {\em Characterizations of finite
  classical polar spaces by intersection numbers with hyperplanes and spaces of
  codimension 2}, Combinatorica, 30 (2010), pp.~25--45.

\bibitem{Graphviz}
{\sc E.~R. Gansner and S.~C. North}, {\em An open graph visualization system
  and its applications to software engineering}, SOFTWARE - PRACTICE AND
  EXPERIENCE, 30 (2000), pp.~1203--1233.

\bibitem{GAP4}
{\sc The GAP~Group}, {\em {GAP -- Groups, Algorithms, and Programming, Version
  4.7.5}}, 2014.

\bibitem{KL90}
{\sc P.~Kleidman and M.~Liebeck}, {\em The subgroup structure of the finite
  classical groups}, vol.~129 of London Mathematical Society Lecture Note
  Series, Cambridge University Press, Cambridge, 1990.

\bibitem{LV13}
{\sc M.~Lavrauw and G.~Van~de Voorde}, {\em Field reduction and linear sets in
  finite projective geometry}, in Topics in Finite Fields, G.~Kyureghyan,
  G.~Mullen, and A.~Pott, eds., vol.~632, AMS Contemporary Math., 2015,
  pp.~271--293.

\bibitem{orb}
{\sc J.~M\"uller, M.~Neunh\"offer, and F.~Noeske}, {\em {orb -- Methods to
  enumerate Orbits, Version 4.7.1}}, 2014.

\bibitem{cvec}
{\sc M.~Neunh\"offer}, {\em {cvec -- Compact vectors over finite fields,
  Version 2.4}}, 2014.

\bibitem{genss}
{\sc M.~Neunh\"offer and F.~Noeske}, {\em {genss -- Generic Schreier-Sims,
  Version 1.6.2}}, 2014.

\bibitem{design}
{\sc L.~Soicher}, {\em {DESIGN -- The Design Package for GAP, Version 1.6}},
  2011.

\bibitem{grape}
\leavevmode\vrule height 2pt depth -1.6pt width 23pt, {\em {GRAPE -- GRaph
  Algorithms using PErmutation groups, Version 4.6.1}}, 2012.

\bibitem{Tits59}
{\sc J.~Tits}, {\em Sur la trialit\'e et certains groupes qui s'en
  d\'eduisent}, Inst. Hautes \'Etudes Sci. Publ. Math.,  (1959), pp.~13--60.

\end{thebibliography}

\end{document}